\documentclass[preprint,12pt]{elsarticle}

\usepackage{amssymb,amsmath}
\usepackage{graphicx}
\usepackage{graphics}
\usepackage{subfigure}
\usepackage{mathrsfs}
\usepackage{eucal}
\usepackage{amsmath}
\usepackage{amssymb}
\usepackage{amsfonts}
\usepackage{multirow}
\usepackage{lscape}

\newcommand{\Eref}[1]{Equation (\ref{#1})}
\newcommand{\fref}[1]{Figure (\ref{#1})}

\newcommand{\uu}{\mathbf{u}}
\newcommand{\xx}{\mathbf{x}}
\newcommand{\KK}{\mathbf{K}}
\newcommand{\BB}{\mathbf{B}}
\newcommand{\DD}{\mathbf{D}}
\newcommand{\qq}{\mathbf{q}}
\newcommand{\aaa}{\mathbf{a}}
\newcommand{\bb}{\mathbf{b}}

\newcommand{\ff}{\mathbf{f}}

\newtheorem{def2}{Definition}

\journal{IJNME}

\begin{document}




\begin{frontmatter}

\title{Integrating strong and weak discontinuities without integration subcells and example applications in an XFEM/GFEM framework} 

\author[a]{Sundararajan~Natarajan}
\author[b]{D Roy Mahapatra}
\author[c]{St\'ephane PA~Bordas\fnref{label2}\corref{cor1}}

\address[a]{Cardiff School of Engineering Theoretical, Applied and Computational Mechanics, Cardiff University, Wales, U.K.}
\address[b]{Assistant Professor, Department of Aerospace Engineering, Indian Institute of Science, Bangalore-560012, INDIA.}
\address[c]{Professor, Cardiff School of Engineering Theoretical, Applied and Computational Mechanics, Cardiff University, Wales, U.K.}

\fntext[4]{Cardiff School of Engineering, Theoretical, Applied and Computational Mechanics, Cardiff University, Wales, U.K. Email: stephane.bordas@alumni.northwestern.edu.
Tel. +44 (0)29 20875941.}

\begin{abstract}
Partition of unity methods, such as the extended finite element
method (XFEM) allow discontinuities to be simulated independently of
the mesh~\cite{Belytschko1999}. This eliminates the need for the
mesh to be aligned with the discontinuity or cumbersome re-meshing,
as the discontinuity evolves. However, to compute the stiffness
matrix of the elements intersected by the discontinuity, a
subdivision of the elements into quadrature subcells aligned with
the discontinuity is commonly adopted. In this paper, we use a
simple integration technique, proposed for polygonal
domains~\cite{Natarajan2009} to suppress the need for element
subdivision. Numerical results presented for a few benchmark
problems in the context of linear elastic fracture mechanics and
a multi-material problem, show that the proposed method
yields accurate results. Owing to its simplicity, the proposed
integration technique can be easily integrated in any existing code.
\end{abstract}


\end{frontmatter}

\section{INTRODUCTION}
The classical finite element method (FEM) is one of the clear
choices to solve problems in engineering and science. But the
classical FEM approaches fail or are computationally expensive for
some classes of problems, viz., equations with rough coefficients
and discontinuities (arising e.g. in the reaction-diffusion
equation, the advection-diffusion equation, crack growth problems,
composites, materials with stiffeners etc.) and problems with
highly oscillatory solutions viz., solution of Helmholtz's equation.
In an effort to improve the FEM, Babu\v{s}ka \textit{et
al.,}~\cite{Babuvska1994} showed that the choice of non-polynomial
ansatz functions, when tailored to the problem
formulation, lead to optimal convergence, whereas the
classical FEM, relying on the approximation properties of
polynomials, performs poorly. Also in the case of
Helmholtz's equation, Melenk~\cite{Melenk:1995} showed that plane
waves displaying the same oscillatory behavior as the solution can
serve as effective enrichment functions. This lead to the birth of
the Partition of Unity Method (PUM).

Belytschko's group in 1999~\cite{Belytschko1999,Moes1999},
exploited the idea of partition of unity enrichment of finite
elements (Babu\v{s}ka \textit{et al.,}~\cite{Babuvska1994}) to solve
linear elastic fracture mechanics problems with minimal remeshing. The
resulting method, known as XFEM is classified as one of the
partition of unity methods in which the main idea is to extend a
classical approximate solution basis by a set of enrichment
functions that carry information about the character of the
solution (e.g., singularity, discontinuity,
boundary layer). 

As it permits arbitrary functions to be incorporated in the
FEM or the mesh-free approximation, partition of unity
enrichment~\cite{Strouboulis2000,Simone2006} leads to greater
flexibility in modeling moving boundary problems, without changing
the underlying mesh 
while the set of
enrichment functions evolve (and/or their supports) with the interface
geometry.


In PU type methods, the enrichment is extrinsic and resolved through additional
degrees of freedom. The enrichment can also be intrinsic, based on the
recent work by Fries and Belytschko~\cite{Fries2006}. In this paper,
we focus on the extrinsic partition of unity enrichment and in general, the field
variables are approximated by~\cite{Belytschko1999,Melenk:1995,Duarte2000,Babuvska1997,Strouboulis2000,Babuvska2003,Simone2006,Hu2007}:

\begin{equation}
\uu^h(\xx)=\sum\limits_{I \in \mathcal{N}^{\rm{fem}}} {N_I
}(\xx)\qq_I + \textup{enrichment functions} \label{eqn:xfemappr}
\end{equation}

\noindent where $N_I(\xx)$ are standard finite element shape
functions, $\qq_I$ are nodal variables associated with node $I$.

XFEM, one of the aforementioned partition of unity methods, was
successfully applied for crack propagation and other fields in
computational physics
~\cite{Chessa2002,Chessa2003,Chopp:2003,Duddu2008,jichoppdolbow,
dolbowmerle,wagner:xfem_stokes,Liu2004,Bos2009,Lecampion2009} and
recently open source XFEM codes were released to help the
development of the method~\cite{bordasc} and numerical
implementation and efficiency aspects were studied~\cite{bordasd}.
XFEM is quite a robust and popular method which is now used for
industrial problems~\cite{Menk2009} and under implementation by
leading computational software companies. It is not the scope of
this paper to review recent advances of partition of unity methods,
and the interested readers are referred to the literature, for
instance, Bordas and Legay~\cite{bordas2006mfm}, Karihaloo
\textit{et al.,}~\cite{Karihaloo2003} and Belytschko \textit{et
al.,}~\cite{Belytschko2009}.

Although XFEM is robust and applied to a wide variety of moving
boundary problems, the flexibility provided by this class of methods
also leads to associated difficulties:
\begin{itemize}
\item when the approximation is discontinuous or non-polynomial in an
element, special care must be taken for numerical integration;
\item the low order of continuity of the solution leads to poor
accuracy (esp. in 3D) of the derivatives close to regions of high
gradient, such as crack fronts~\cite{Xiao5} which motivated recent work on adaptivity for GFEM~\cite{Barros2004,Strouboulis2006}, meshfree methods~\cite{Belytschko9,Rabczuk1,Rabczuk3}
and XFEM (Bordas and Duflot~\cite{xmls1}, Bordas {\it et
al.}~\cite{xmls0}, R\'odenas {\it et al.} \cite{Rodenas2008}).
\end{itemize}

An important first attempt to 
simplify numerical integration was by Ventura~\cite{Ventura2006},
who focuses on the elimination of quadrature subcells commonly
employed to integrate strongly or weakly discontinuous and
non-polynomial functions present in the enriched FE approximation.
His work is based on replacing non-polynomial functions by
`equivalent' polynomials. The proposed method is exact for
triangular and tetrahedral elements, but for quadrilateral elements,
when the opposite sides are not parallel, additional approximation
is introduced.

Another method that alleviates this difficulty is strain smoothing~\cite{Bordas2009,Bordas2008}. The main idea is to combine
the smoothed finite element method
(SFEM)~\cite{liu2007,Nguyen-Xuan2008} with the XFEM. The SFEM relies
on strain smoothing, which was proposed by Chen \textit{et
al.,}~\cite{Chen2001} for meshless methods, where the strain is
written as the divergence of a spatial average of the standard
(compatible) strain field -- i.e., symmetric gradient of the
displacement field. In strain smoothing, the surface
integration is transformed into an equivalent boundary integration
by use of the Green-Ostrogradsky theorem. The Smoothed XFEM 
was introduced in~\cite{Bordas2009}, but much remains to be understood regarding the convergence, stability and accuracy of this method.

 In this paper, we propose to use the new numerical integration technique
proposed by the authors for arbitrary
polygonal domains~\cite{Natarajan2009,Natarajan2009a,Natarajan2008} to
compute the stiffness matrix. Each part of the elements that are cut
or intersected by a discontinuity is conformally mapped onto a unit
disk using Schwarz-Christoffel mapping. A midpoint quadrature is
used to obtain the integration points as opposed to the regular
Gau\ss~cubature rule. Thus, the proposed method which works only in 2D, eliminates
the need to sub-divide the elements cut by discontinuities into
quadrature subcells for the purpose of numerical integration of the
stiffness matrix. The Schwarz-Christoffel mapping has been applied to mesh free
Galerkin method by Balachandran \textit{et al.,}~\cite{Balachandran2008} to obtain the
weight function of the arbitrary shaped support domain obtained from natural neighbor algorithm.

The paper is organized as follows. In the next section, we briefly
recall the basic equations of the XFEM. Section~\ref{scmapp} will
explain the new numerical integration scheme. The efficiency and
convergence properties of the proposed method are illustrated in
section~\ref{numeexamples} with a few benchmark problems taken from
linear elastic fracture mechanics and a multi-material problem,
which is followed by some concluding remarks in the last section.

\section{BASICS OF PARTITION OF UNITY METHODS FOR DISCONTINUITIES AND SINGULARITIES}
\label{xfem:basics} In this section, we give a brief overview of
partition of unity methods in FEM for problems with strong and weak
discontinuities
We focus on the extended
FEM~\cite{Belytschko1999,bordas2006mfm}, but the method is
identical for alternatives such as the
GFEM~\cite{Barros2004,Strouboulis2006}.

\begin{figure}[htpb]
\centering
\scalebox{0.7}{\input{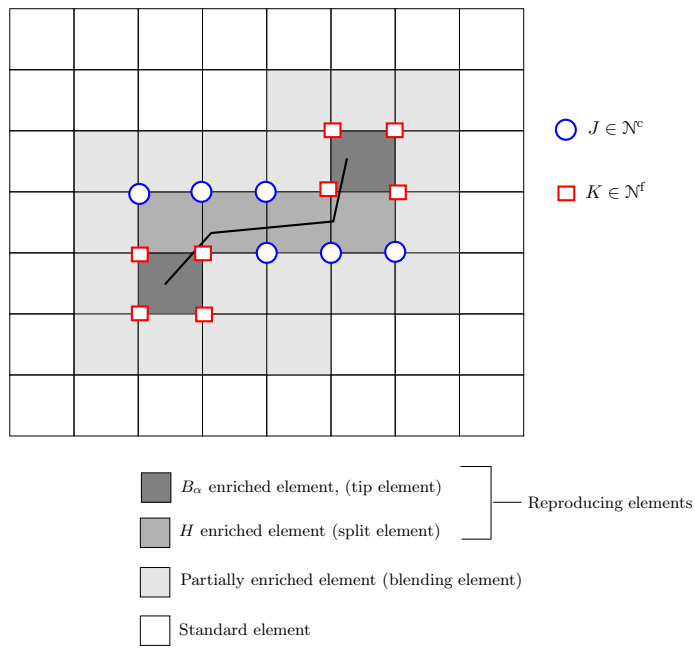}}
\caption{A typical FE mesh with a discontinuity. The squared nodes are enriched with the Heaviside function, and the circled nodes are enriched with near-tip asymptotic field obtained from Westergaard solution~\cite{Belytschko1999}} \label{fig:xfem_crack}
\end{figure}

With a regular finite element method, the mesh has to
conform to the discontinuities and a very fine mesh is required in
regions with sharp gradients. When the discontinuity surface
evolves, cumbersome remeshing is required. The XFEM alleviates these
difficulties by allowing the discontinuities to be independent of
the mesh. An XFEM model consists of a regular FE mesh, which is
independent of the discontinuity geometry. \fref{fig:xfem_crack}
illustrates a typical FE mesh with a crack.

The main idea is to extend the approximation basis by a set of
enrichment functions, that are chosen based on the local behavior of
the problem. For the case of linear elastic fracture mechanics, two
sets of functions are used: a Heaviside jump function to capture the
jump across the crack faces and asymptotic branch functions that
span the two-dimensional asymptotic crack tip fields. The enriched
approximation for fracture mechanics problems takes the
form~\cite{Belytschko1999,bordas2006mfm,xfemCourseNotes}:

\begin{equation}
\uu^h(\xx)=\sum\limits_{I \in \mathcal{N}^{\rm{fem}}} {N_I
}(\xx)\qq_I + \sum\limits_{J \in
\mathcal{N}^{c}}{{N_J}}(\xx)H(\xx)\aaa_J + \sum\limits_{K \in
\mathcal{N}^{f}}{{N_K}} (\xx)\sum_{\alpha=1}^4
B_{\alpha}(\xx)\bb^{\alpha}_K ,\label{eqn:uS2}
\end{equation}

\noindent where $\aaa_J$ and $\bb_K$ are nodal degrees of freedom
corresponding to the Heaviside function $H$ and the near-tip
functions, $\{B_\alpha\}_{1 \leq \alpha \leq 4}$, given by:

\begin{equation}
\{B_\alpha\}(r,\theta)_{1 \leq \alpha \leq 4} = \sqrt{r} \left\{
\sin\left(\frac{\theta}{2}\right),
\cos\left(\frac{\theta}{2}\right), \sin\left(\theta\right)
\sin\left(\frac{\theta}{2}\right), \sin\left(\theta\right)
\cos\left(\frac{\theta}{2}\right) \right\}.
\end{equation}

\noindent Nodes in set $\mathcal{N}^c$ are such that their support
is split by the crack and nodes in set $\mathcal{N}^f$ belong to the
elements that contain a crack tip. These nodes are enriched with the
Heaviside and near-tip (branch functions) fields, respectively. In
the discretization of~\Eref{eqn:uS2}, the displacement field is
global, but the supports of the enriching functions are local
because they are multiplied by the nodal shape functions.

This modification of the displacement approximation does not
introduce a new form of the discretized finite element equilibrium
equation, but leads to an enlarged problem to solve:

\begin{equation}
\renewcommand\arraystretch{2}
\left[\begin{array}{*{20}c}
                \KK_{uu}  & \KK_{ua}  & \KK_{ub}\\
                 \KK_{au}  &  \KK_{aa} & \KK_{ab} \\
                 \KK_{bu} & \KK_{ba} & \KK_{bb}
                 \end{array} \right] \left\{ \begin{array} {c} \qq \\ \aaa \\ \bb \end{array} \right\}
= \left\{ \begin{array} {c} \ff_q \\ \ff_a  \\ \ff_b\end{array} \right\}
\end{equation}

\noindent where the element stiffness matrix is given by:

\begin{equation}
\renewcommand\arraystretch{2}
\left[\begin{array}{*{20}c}
                \KK_{uu}^e  & \KK_{ua}^e & \KK_{ub}^e \\
                 \KK_{au}^e  &  \KK_{aa}^e & \KK_{ab}^e \\
                 \KK_{bu}^e  &  \KK_{ba}^e & \KK_{bb}^e
                 \end{array} \right] =  \int_{\Omega^e} \left[ \begin{array}{*{20}c}
                \BB_{std}^T \DD \BB_{std} &  \BB_{std}^T \DD
                 \BB_{enr}^1 &  \BB_{std}^T \DD
                 \BB_{enr}^2\\  (\BB_{enr}^1)^T \DD
                 \BB_{std} &  (\BB_{enr}^1)^T \DD
                 \BB_{enr}^1  &  (\BB_{std}^1)^T \DD
                 \BB_{enr}^2 \\
                  (\BB_{enr}^2)^T \DD
                 \BB_{std} &  (\BB_{enr}^2)^T \DD
                 \BB_{enr}^1  &  (\BB_{std}^2)^T \DD
                 \BB_{enr}^2 \end{array} \right]~d\Omega_e
\end{equation}

\noindent where $\BB_{std}$ is the standard strain-displacement matrix, $\BB_{enr}^1$ and $\BB_{enr}^2$  are the enriched parts of the strain-displacement matrix corresponding to the Heaviside and asymptotic functions, respectively. $\DD$ is the
material matrix.

\begin{figure}[htpb]
\centering
\scalebox{0.7}{\input{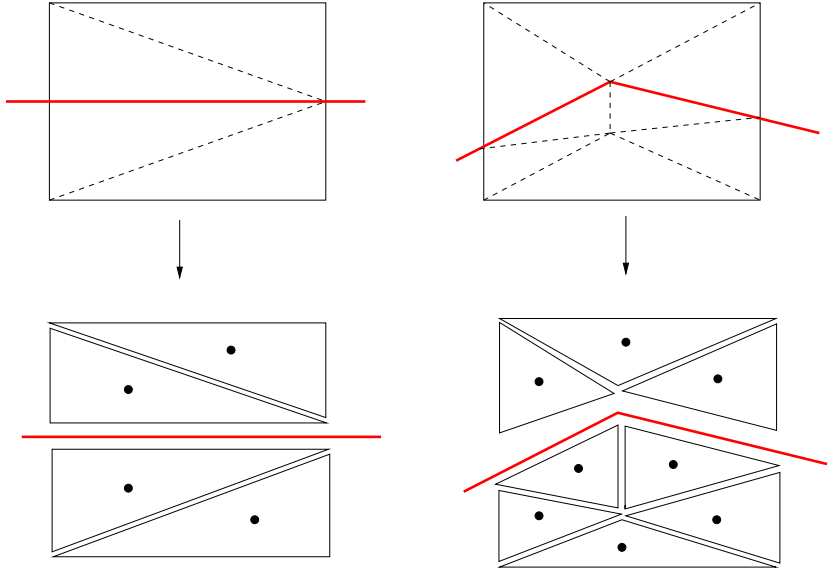}}
\caption{Integration in an element with a straight and kinked
discontinuity: standard decomposition of an element for integration
of a discontinuous weak form for XFEM. Gau\ss~points are introduced
within each triangle to ensure proper integration of the
discontinuous displacement field.} \label{fig:xfem_idea}
\end{figure}

The numerical integration of the stiffness matrix in elements
intersected by a discontinuity, be it a material interface (weak
discontinuity) or a crack (strong discontinuity) is not trivial. The
standard Gau\ss~quadrature cannot be applied in elements enriched by
discontinuous terms, because the Gau\ss~quadrature implicity assumes
a polynomial approximation. This is circumvented by partitioning the
elements into subcells aligned to the discontinuity surface, in
which the integrands are continuous and differentiable
(see~\fref{fig:xfem_idea}). Although, the generation of quadrature
subcells does not alter the approximation properties, it inherently
introduces a `mesh' requirement. The steps involved in this
approach are: (1) Split the element into subcells with the subcells
aligned to the discontinuity surface, usually the subcells are
triangular (see~\fref{fig:xfem_idea}) and (2) numerical integration
is performed with the integration points from triangular quadrature.
The subcells must be aligned to the crack or interface and this is
costly and less accurate if the discontinuity is curved. Similar attempts were made to improve the integration of discontinuities in meshfree methods~\cite{Chen2,Rabczuk10,Rabczuk6,Rabczuk13,Bordas1}. To
alleviate this difficulty, we propose to use the new numerical
integration technique proposed by the authors~\cite{Natarajan2009}
for polygonal domains to integrate over the elements intersected by
the discontinuity. The next section will briefly review the
Schwarz-Christoffel conformal mapping (SCCM) and then discuss how to
numerically integrate over elements intersected by discontinuities
using the SCCM.


\section{SCHWARZ-CHRISTOFFEL CONFORMAL MAPPING}
\label{scmapp} Conformal mapping is extremely important in complex
analysis and finds its application in many areas of physics and
engineering. A conformal transformation or biholomorphic map is a
transformation that preserves local angles. In other words, if
$\Gamma_1$ and $\Gamma_2$ are two curves that intersect at an angle
$\theta_z$ in the $z-$ plane at a point $p$, then the images
$f\left(\Gamma_1\right)$ and $f\left(\Gamma_2\right)$ intersect at
an angle $\theta_w = \theta_z$ at $q=f\left(p\right)$. A
Schwarz-Christoffel mapping is a transformation of the complex plane
that maps the upper half-plane conformally to a polygon.

\begin{def2} A Schwarz-Christoffel map is a function $f$ of
the complex variable that conformally maps a canonical domain in the
$z-$ plane (a half-plane, unit disk, rectangle, infinite strip) to a
`closed' polygon in the $w-$ plane.
\end{def2}

Consider a polygon in the complex plane. The Riemann mapping theorem
implies that there exists a bijective holomorphic mapping $f$ from
the upper half plane $\left\{ \zeta \in \mathbb{C}:
\operatorname{Im}\,\zeta > 0 \right\}$ to the  interior of the
polygon.  The function $f$ maps the real axis to the edges of the
polygon. If the polygon has interior angles
$\alpha,\beta,\gamma,\ldots$, then this mapping is given by,

\begin{equation}
f(\zeta) = \int^\zeta
\frac{K}{(w-a)^{1-(\alpha/\pi)}(w-b)^{1-(\beta/\pi)}(w-c)^{1-(\gamma/\pi)}
\cdots} \,\mbox{d}w \end{equation}

\noindent where $K$ is a constant, and $a < b < c < ...$ are the
values, along the real axis of the $\zeta$ plane, of points
corresponding to the vertices of the polygon in the $z$ plane. A
transformation of this form is called a Schwarz-Christoffel mapping.

\begin{figure}[h]
\centering
\includegraphics[angle=0,width=1\textwidth]{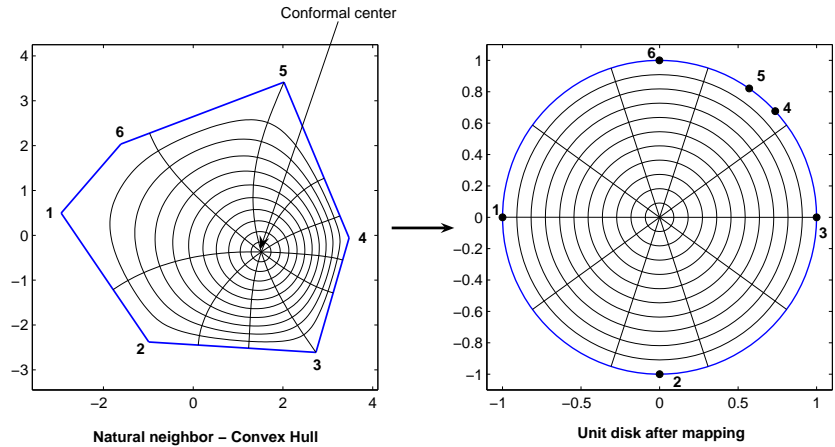}
\caption{Mapping of the physical domain to the unit disk. This
figure was produced with the MATLAB SC Toolbox~\cite{Driscoll1996}}
\label{fig:conformal}
\end{figure}

\begin{figure}[htpb]
\centering
\subfigure[Midpoint rule]{\input{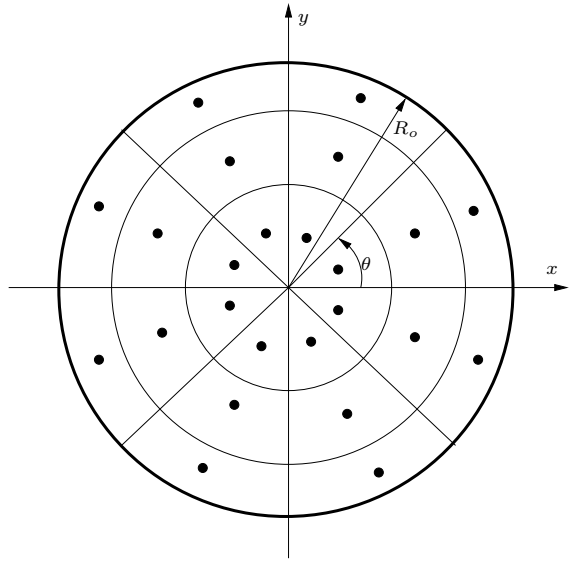}}
\subfigure[Gau\ss{}-Chebyshev rule]{\input{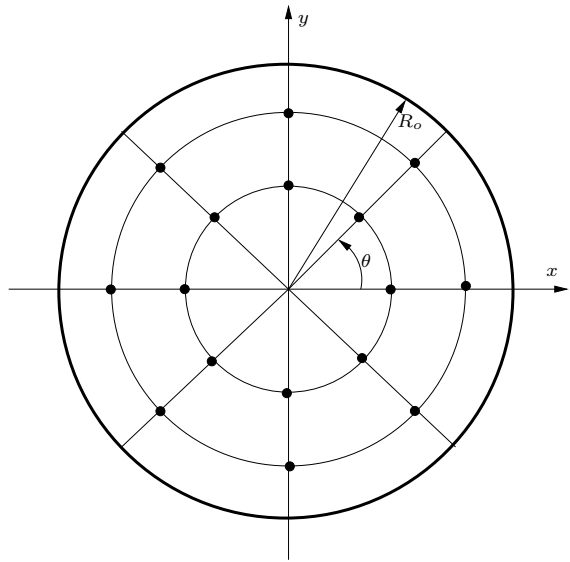}}
\caption{Quadrature rules on a disk} \label{fig:points}
\end{figure}

Recently, the authors proposed a new numerical integration
technique~\cite{Natarajan2009,Natarajan2009a,Natarajan2008} using
the Schwarz-Christoffel mapping and cubature rules on a disk to
numerically integrate over arbitrary polygons that arise in
polygonal finite element methods. \fref{fig:conformal} shows
the conformal mapping of an arbitrary polygon onto a unit disk on
which a midpoint rule~\cite{Suvrau2001} or Gau\ss~Chebyshev
rule~\cite{Peirce1957} is used to obtain integration points. The
distributions of integration points of the mid point quadrature and
Gau\ss{}-Chebyshev quadrature are illustrated in \fref{fig:points}.

\begin{figure}[htpb]
\centering
\input{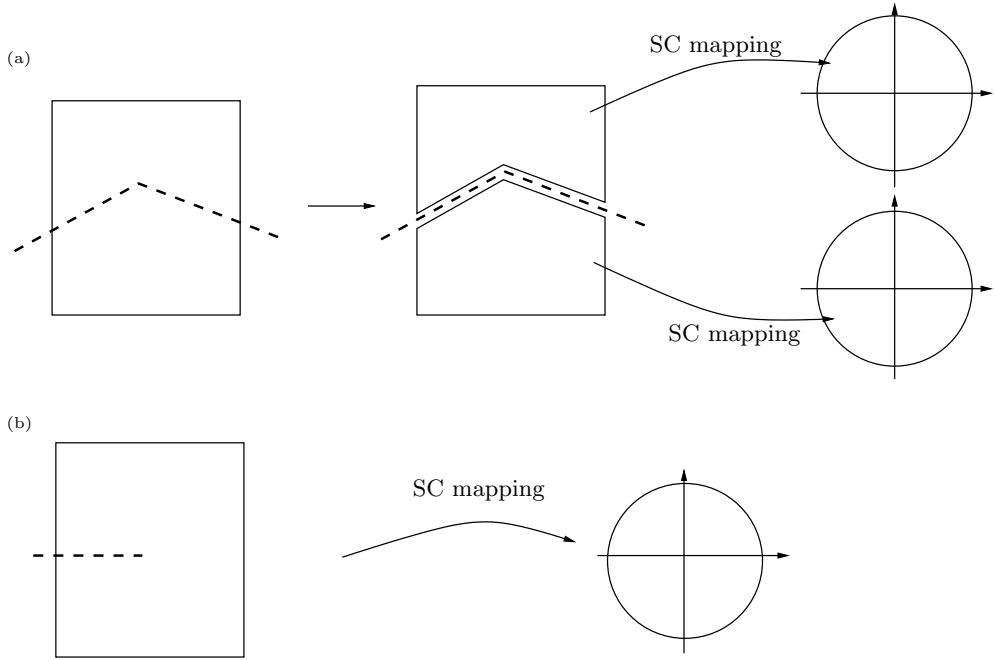}
\caption{Integration over an element with discontinuity (dotted line):
(a) with kinked discontinuity, representing the split element and
(b) strong discontinuity, representing the tip element. In both cases, the sub-polygon is mapped
conformally onto the unit disk using Schwarz-Christoffel conformal
mapping.} \label{fig:split}
\end{figure}

In this paper, we propose to use conformal mapping in the
context of the XFEM for 2D problems to numerically integrate over
elements where the approximation or its derivatives is
discontinuous. Each part of the element is conformally mapped onto a
unit disk using the technique proposed in~\cite{Natarajan2009}.
\fref{fig:split} illustrates the above idea for split and tip
elements.

One of the cubature rules mentioned above is used to obtain
integration points. For elements that are not intersected by a
discontinuity surface, the standard isoparametric mapping is
implemented with four integration points for each element.

\clearpage

\section{NUMERICAL EXAMPLES}
\label{numeexamples}

In this section, we illustrate the effectiveness of the
proposed method by solving a few benchmark problems taken from
linear elastic fracture mechanics and a multi-material problem. We
first consider an infinite plate under tension, then a plate with an
edge crack, for which analytical solutions are available.
Then, a multiple crack problem, followed by a
multi-material problem is considered  and as a last example crack growth in a double cantilever beam is studied. In this study, unless
otherwise mentioned, quadrilateral elements\footnote{Bilinear
element, 4 noded quadrilateral element} are used. In case of the
XFEM with a standard integration approach, 13 integration points per
subcell are used and a similar number of integration points are used
for the SCCM, i.e., for a tip element with six subcells, 78 (=13 $\times
6$) integration points are used for both methods. In this study, we use only topological enrichment, i.e., only the tip element is enriched by near tip functions~\cite{B'echet2005,Laborde2005,Chahine2008}.

\subsection{Infinite plate under tension}
\begin{figure}[htbp]
\centering
\scalebox{0.6}{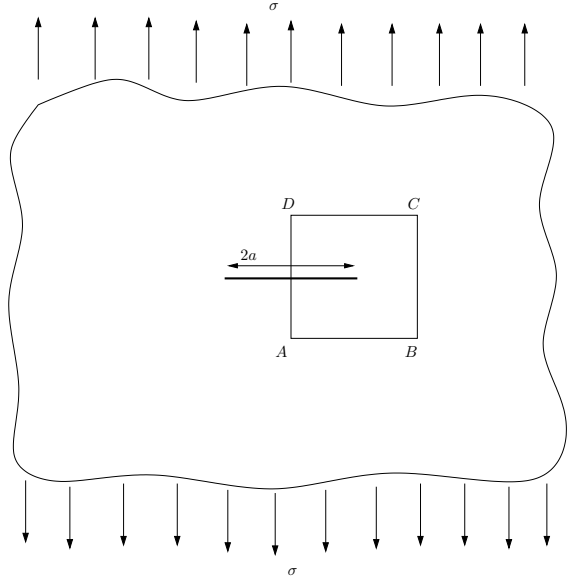}
\caption[Infinite cracked plate]{Infinite cracked plate under remote
tension: geometry and loads } \label{fig:infPlate}
\end{figure}

Consider an infinite plate containing a straight crack of length $a$
and loaded by a remote uniform stress field $\sigma$ as shown
in~\fref{fig:infPlate}. Along ABCD the closed form solution in terms
of polar coordinates in a reference frame $(r,\theta)$ centered at
the crack tip is

\begin{subequations}
\begin{eqnarray}
\label{eqn:infinite} \sigma_x (r,\theta) = \frac{K_I}{\sqrt{r}}
\cos\frac{\theta}{2} \left(1-\sin\frac{\theta}{2}
\sin\frac{3\theta}{2}\right)
\\
\sigma_y (r,\theta)  = \frac{K_I}{\sqrt{r}} \cos\frac{\theta}{2}
\left(1+\sin\frac{\theta}{2} \sin\frac{3\theta}{2}\right)\\
\sigma_{xy} (r,\theta)  = \frac{K_I}{\sqrt{r}} \sin\frac{\theta}{2}
\cos\frac{\theta}{2} \cos\frac{3 \theta}{2}
\end{eqnarray}
\end{subequations}

The closed form near-tip displacement field is:
\begin{subequations}
\begin{eqnarray}
u_x  (r,\theta) = \frac{2(1+\nu)}{\sqrt{2\pi}} \frac{K_I}{E}
\sqrt{r} \cos\frac{\theta}{2} \left(2-2\nu-\cos^2\frac{\theta}{2}\right) \\
u_y  (r,\theta) = \frac{2(1+\nu)}{\sqrt{2\pi}} \frac{K_I}{E}
\sqrt{r} \sin\frac{\theta}{2}
\left(2-2\nu-\cos^2\frac{\theta}{2}\right)
\end{eqnarray}
\label{eqn:disp}
\end{subequations}

In the two previous expression $K_I=\sigma \sqrt{\pi a}$ denotes the
stress intensity factor (SIF), $\nu$ is Poisson's ratio and $E$ is
Young's modulus. All simulations are performed with $a=100$mm and
$\sigma=10^4$ N/mm$^2$ on a square mesh with sides of length $10$mm.

\begin{figure}[htpb]
\centering
\includegraphics[width=0.7\textwidth]{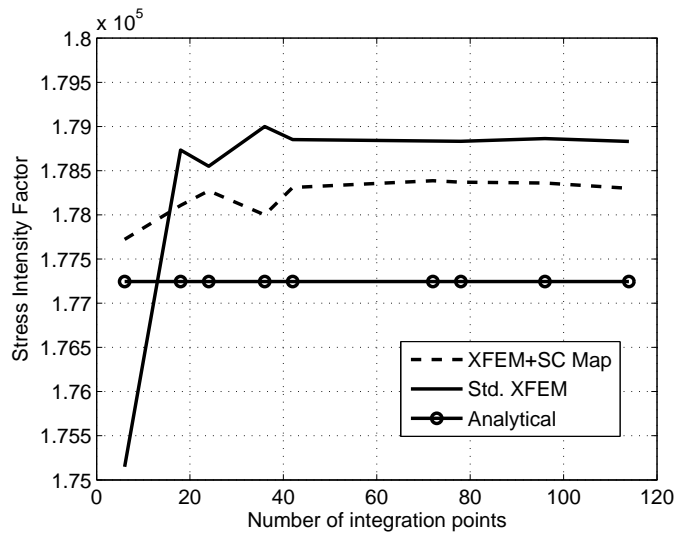}
\caption{Griffith Problem: the convergence of the numerical stress intensity factor with number of integration points in the tip element. A structured quadrilateral mesh $(60\times 60)$ is used.}
\label{fig:sif_conve_ng}
\end{figure}

\begin{figure}[htpb]
\centering
\subfigure[]{\includegraphics[width=0.7\textwidth]{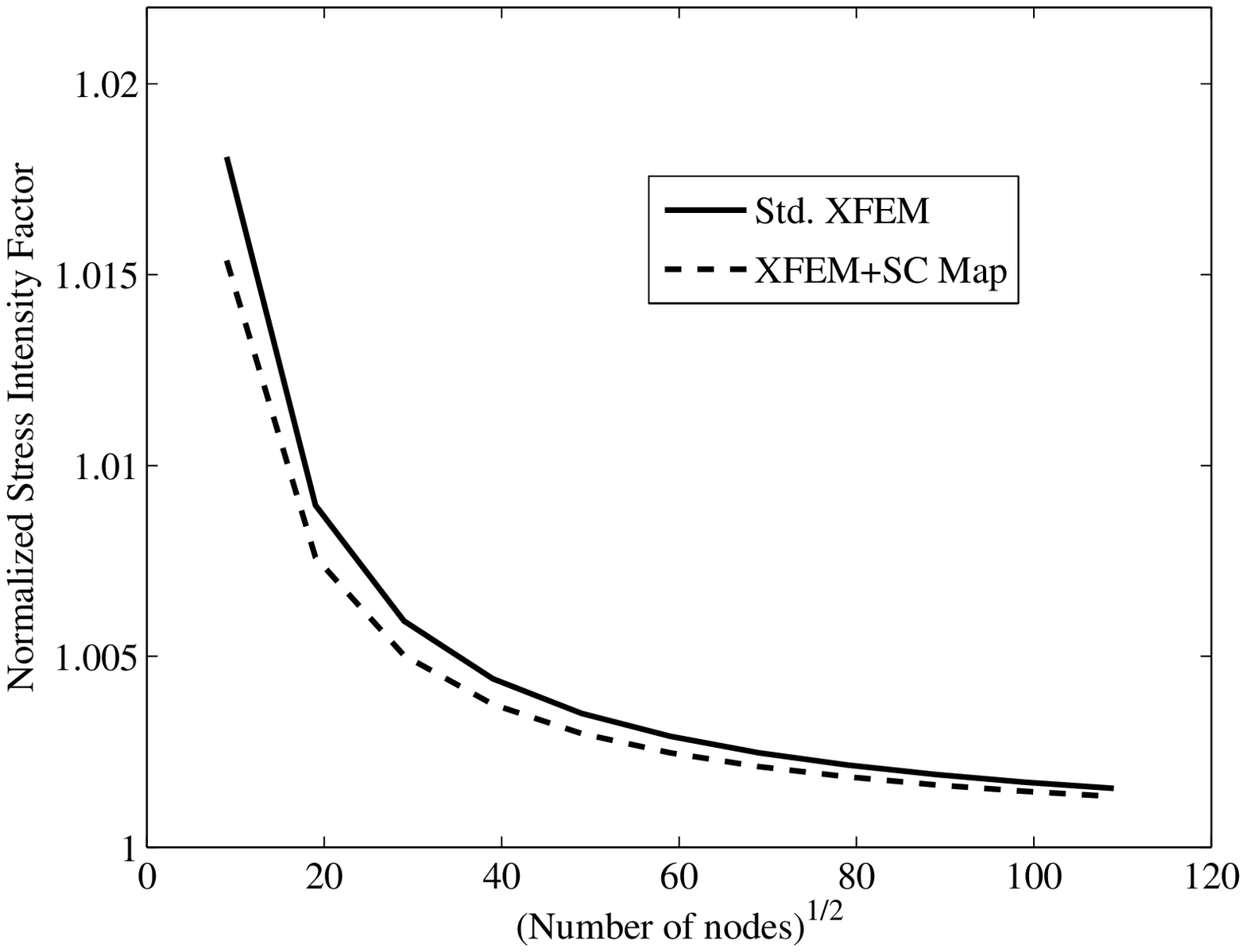}}
\subfigure[]{\includegraphics[width=0.7\textwidth]{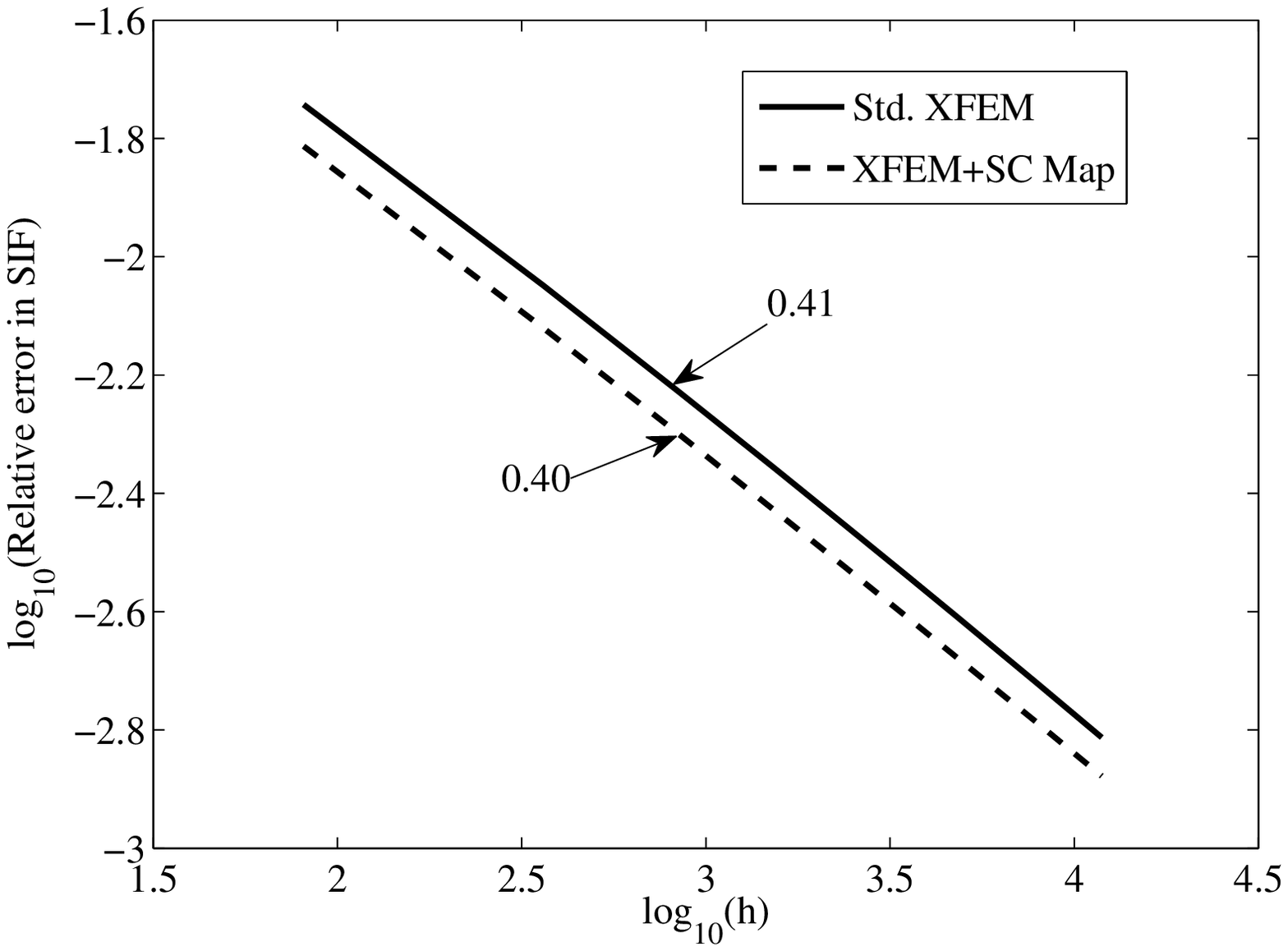}}
\caption{Griffith Problem: the convergence of the numerical stress
intensity factor to the analytical stress intensity factor and
convergence rate.} \label{fig:sif_conve}
\end{figure}

Before carrying out the mesh convergence study and other numerical studies, the influence of the number of integration points in the tip element on the numerical SIF is studied. A structured quadrilateral mesh $(60 \times 60)$ is used for the study and the number of integration points are varied until the difference between two consecutive computations are less than a specified tolerance. During this study, the number of integration points for both methods are kept the same. The convergence of the SIF with the integration points is shown in~\fref{fig:sif_conve_ng}. It is seen that with the increase in the number of integration points, the SIF initially increases but reaches a constant value beyond $60$ integration points. 

The convergence and rate of convergence in numerical stress
intensity factor are shown in~\fref{fig:sif_conve}. It is seen that
for the same number of integration points, the proposed method
outperforms (although only slightly) the conventional XFEM.
But with increase in mesh size, both techniques of numerical
integration approach the analytical solution.

Examining \fref{fig:sif_conve} shows that the convergence rates in the SIFs are suboptimal both for the XFEM with standard integration and the XFEM with the new integration technique proposed. This is for two reasons: (i) only the tip element is enriched (topological enrichment~\cite{Xiao5,B'echet2005,Laborde2005}, which asymptotically reduces the XFEM approximation space to the standard FEM approximation space, this limits the optimal convergence rate to 0.5 in the presence of a square root singularity; (ii)no blending correction \cite{chessablending} is performed, which leads to yet a smaller convergence rate of 0.4, which is consistent with the literature~\cite{Laborde2005,Stazi2003}.

\begin{figure}[htpb]
\centering
\subfigure{\includegraphics[width=0.49\textwidth]{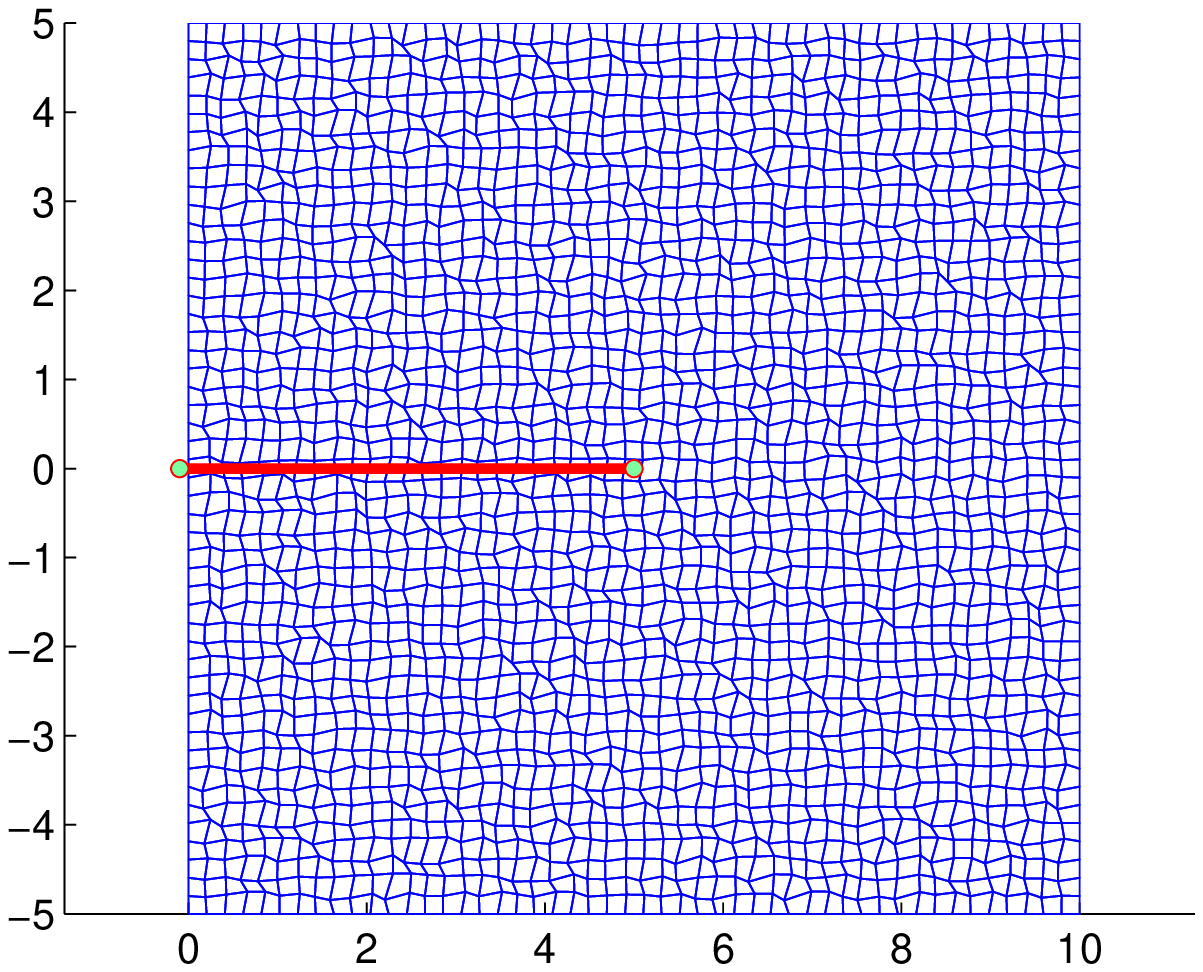}}
\subfigure{\includegraphics[width=0.49\textwidth]{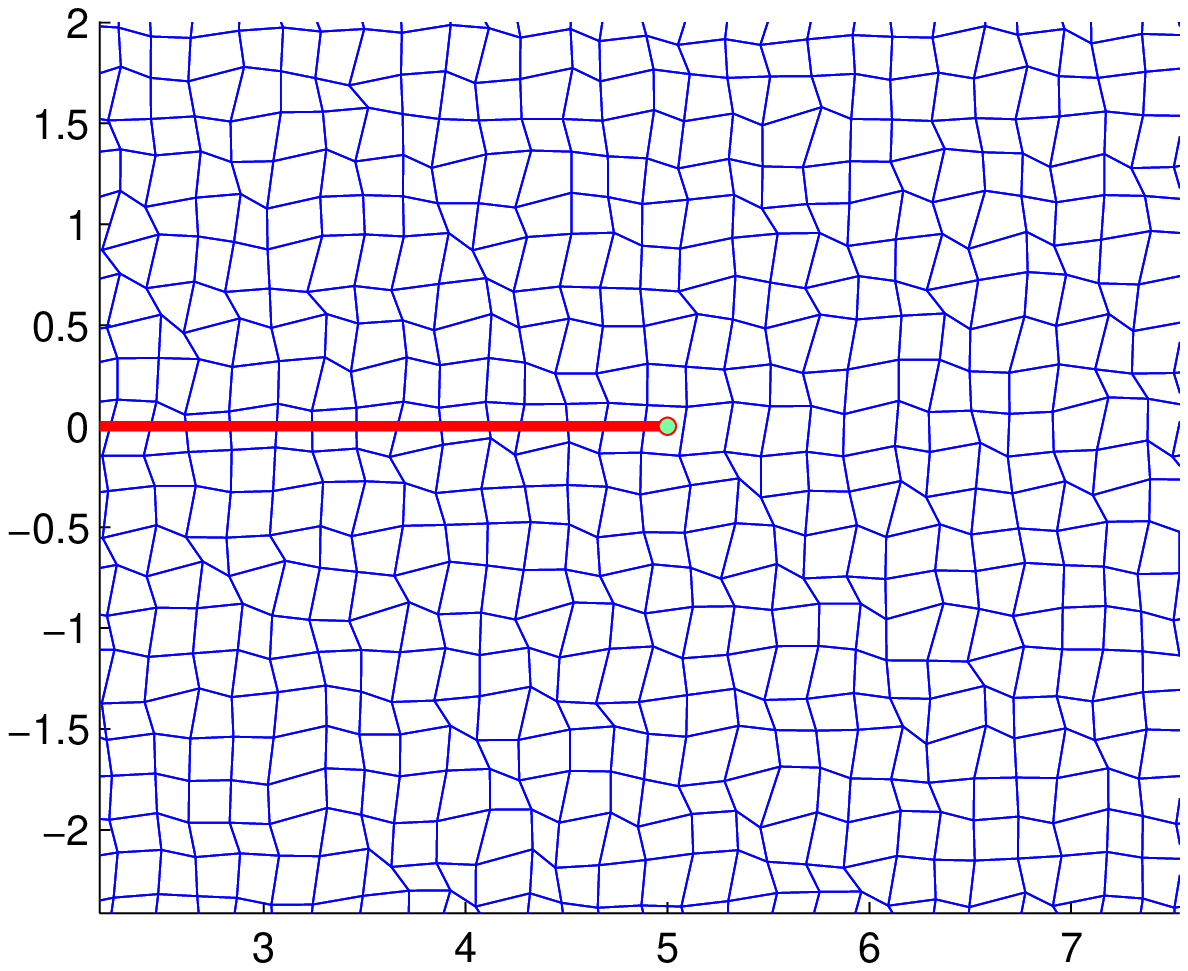}}
\caption{Polygonal mesh for infinite plate problem with a crack under uniform far field tension. The irregular mesh is generated with an irregularity factor, $\alpha_{ir}=0.4$.}
\label{fig:irremesh}
\end{figure}

To demonstrate the effectiveness of the proposed method in case of skewed elements, the domain is meshed with irregular elements. The coordinates of interior nodes are given by

\begin{subequations}
\begin{eqnarray}
x^\prime = x + (2r_c-1) \alpha_{ir} \Delta x \\
y^\prime = y + (2r_c-1) \alpha_{ir} \Delta y
\end{eqnarray}
\end{subequations}

\noindent where $r_c$ is a random number between 0 and 1.0, $\alpha_{ir}$ is an irregularity factor controlling the shapes of the distorted elements and $\Delta x, \Delta y$ are initial regular element sizes in the $x-$ and $y-$directions respectively. The discretization is shown in~\fref{fig:irremesh}.

\begin{figure}[htpb]
\centering
\includegraphics[width=0.7\textwidth]{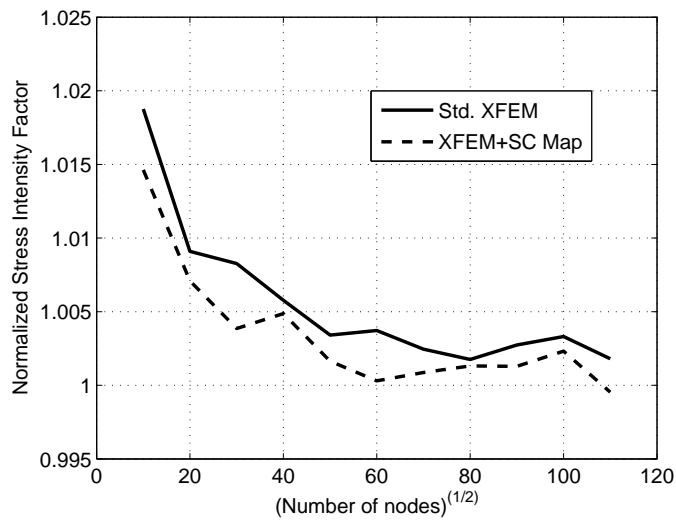}
\caption{Griffith Problem: the convergence of the numerical stress
intensity factor to the analytical stress intensity factor for distorted elements.}
\label{fig:sif_conve_irre}
\end{figure}

The convergence of the numerical stress intensity factor for the irregular mesh is shown in~\fref{fig:sif_conve_irre}. However, the convergence for distorted elements exhibits a non-uniform behavior. Results show that the results from both the numerical techniques are comparable. The advantage of the proposed method is that it eliminates the need to subdivide the split and the tip elements.

\subsection{Edge crack under tension}
\begin{figure}[htp]
\centering
\input{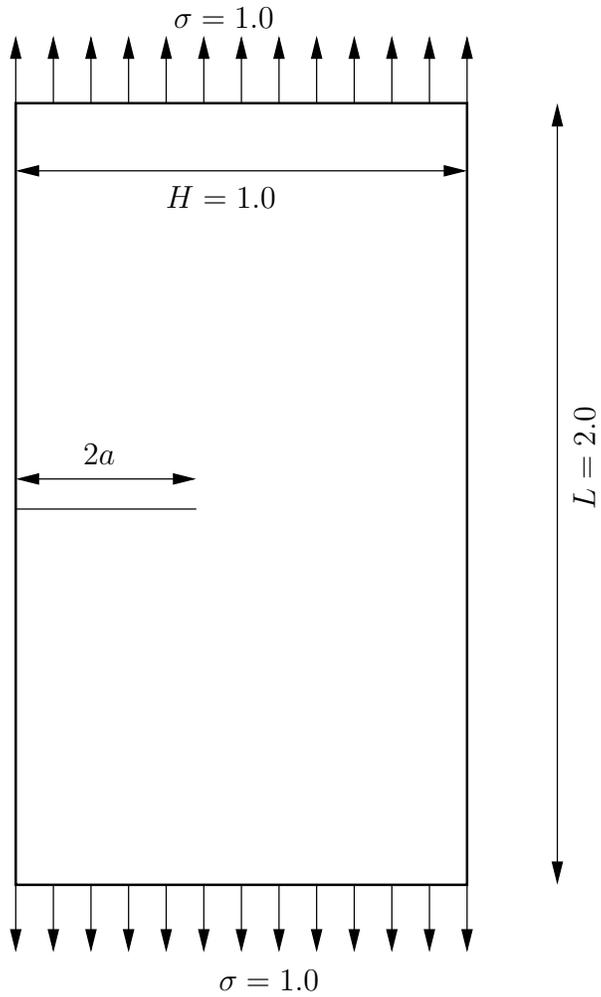}
\caption{Plate with edge crack under tension} \label{fig:edgecrack}
\end{figure}

A plate of dimension 1$\times$2 is loaded by a tension $\sigma=1$
over the top edge. The displacement along the $y$-axis is fixed at
the bottom right corner and the plate is clamped at the bottom left
corner. The geometry, loading, boundary conditions and domain
discretization are shown in Figure~\ref{fig:edgecrack}. The
reference mode I SIF is given by

\begin{equation}
\label{eqn:modeIsif} K_I=F \left( \frac{a}{H} \right) \sigma \sqrt{\pi a}
\end{equation}

\noindent where $a$ is the crack length, $H$ is the plate width and
$F(\frac{a}{H})$ is an empirical function given as (For
$(\frac{a}{H}) \le$ 0.6)

\begin{equation}
\label{eqn:empirical} F\left(\frac{a}{H}\right) = 1.12 -
0.231\left(\frac{a}{H}\right) + 10.55 \left(\frac{a}{H} \right)^2 -
21.72\left(\frac{a}{H}\right)^3 + 30.39\left(\frac{a}{H}\right)^4
\end{equation}

\begin{figure}[htpb]
\centering
\subfigure[]{\includegraphics[width=0.7\textwidth]{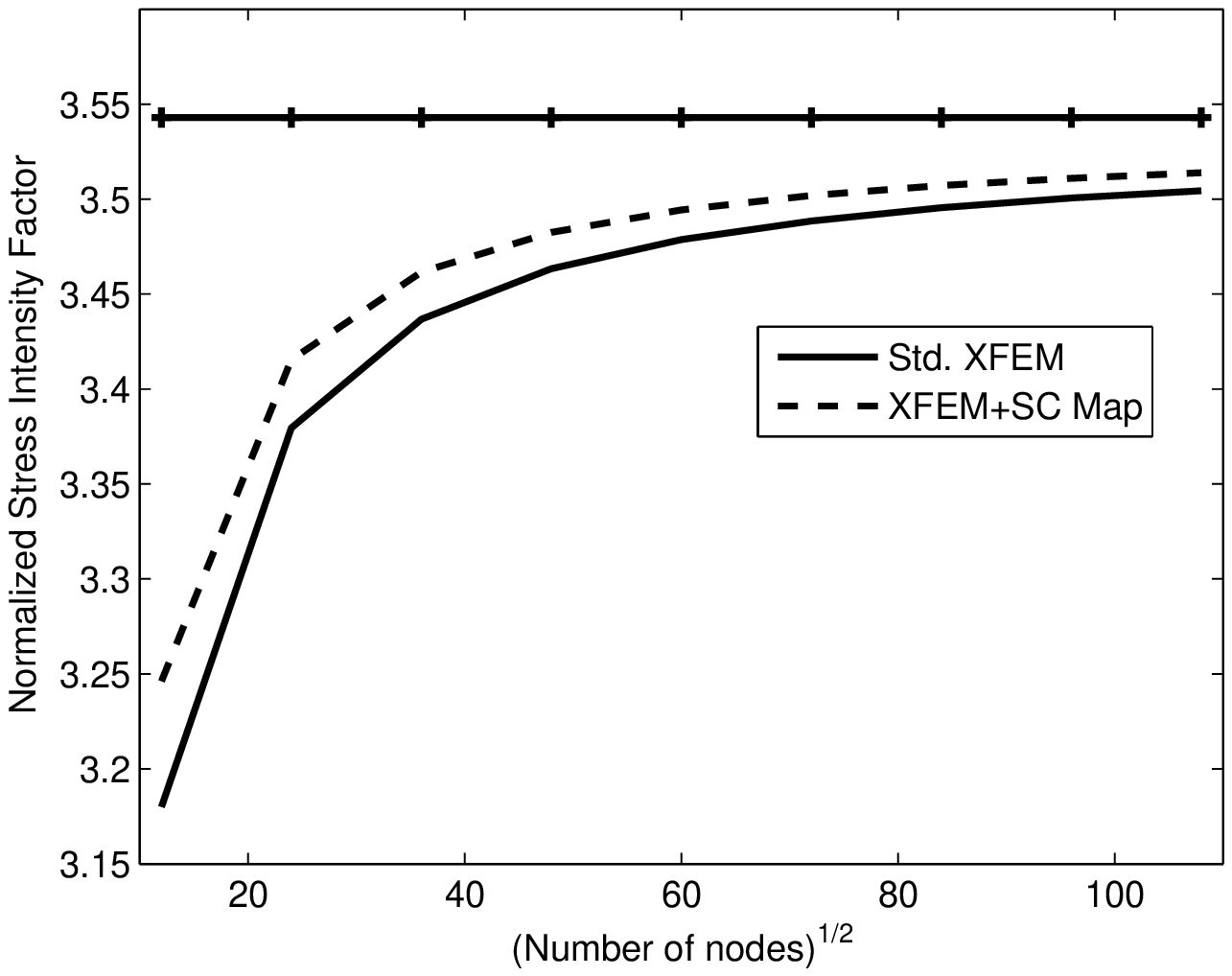}}
\subfigure[]{\includegraphics[width=0.7\textwidth]{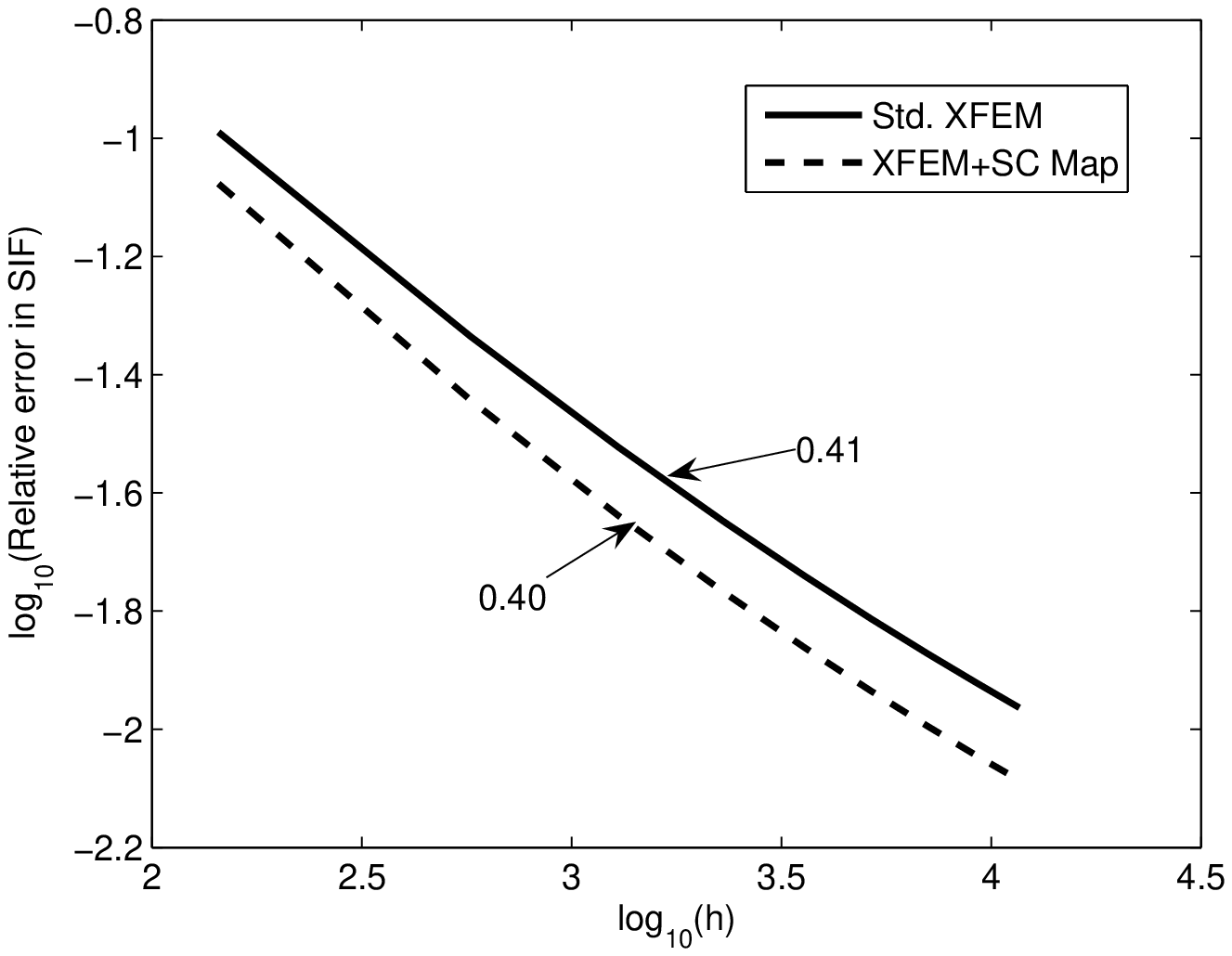}}
\caption{Edge crack problem: the convergence of the numerical stress
intensity factor to the analytical stress intensity factor and
convergence rate. Both the methods show a rate of convergence of 0.4, very close to the ones obtained in~\cite{B'echet2005}.} \label{fig:edg_sif_conve}
\end{figure}

The convergence of the mode I SIF with mesh size and the rate of
convergence of the SIF for a plate with an edge crack is shown
in~\fref{fig:edg_sif_conve}. It is seen that with decrease in mesh
density, both methods approach the analytical solution. Also, the
proposed method performs slightly better than the standard XFEM. The rate of convergence for both methods is 0.4, very similar to the results available in the literature~\cite{Belytschko1999,Stazi2003,B'echet2005}. The rate of convergence can be improved by using a `geometrical' enrichment as suggested by B\'echet \textit{et al.,}~\cite{B'echet2005} or by modifying the enrichment functions such that they are zero in the standard elements and vary continuously in the blending elements as suggested by Fries~\cite{Fries2008}.

In all the above examples, the background FE mesh is made up of regular quadrilateral elements. The numerical integration of the stiffness matrix over regular quadrilateral elements is not a difficult task. The real challenge is when the background mesh is made up of polygons or when the crack faces are irregular, i.e., when the crack faces cut the elements in such a way that at least one of the subdomains, created by the intersection of the geometry with the mesh, is a polygonal with more than 3 edges. To demonstrate the usefulness of the proposed integration technique, we consider the problem of inclined crack in tension with bilinear quadrilateral element as background FE mesh. The study of problems involving polygonal FE meshes and arbitrary crack faces will be the topic of future papers.

\subsection{Inclined crack in tension}
\begin{figure}[htbp]
\centering
\includegraphics[angle=0,width=0.4\textwidth]{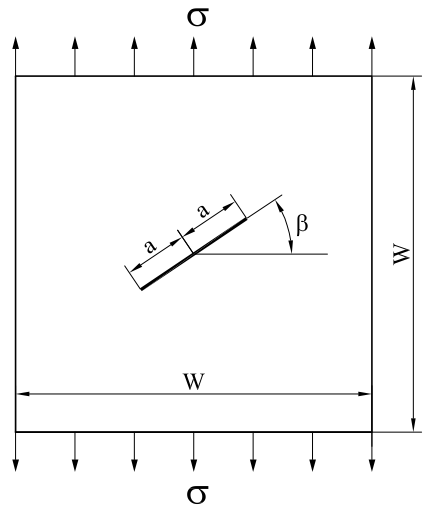}
\caption[]{Inclined crack in tension} \label{fig:incl_crk}
\end{figure}

Consider a plate with an angled crack subjected to a far field
uniaxial stress field (see ~\fref{fig:incl_crk}). In this example, $K_I$ and $K_{II}$ are
obtained as a function of the crack angle $\beta$. For the loads
shown, the analytical stress intensity factors are given by
~\cite{Sih1973}
\begin{eqnarray}
\label{eqn:inc_crk} K_I = \sigma \sqrt{\pi a} \cos\beta \cos\beta,
\quad K_{II} = \sigma \sqrt{\pi a} \cos\beta \sin\beta.
\end{eqnarray}

\begin{figure}[htpb]
\centering
\includegraphics[angle=0,width=0.7\textwidth]{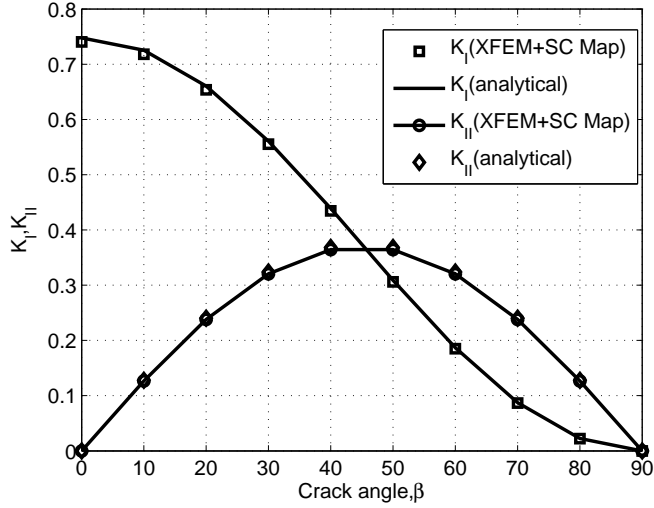}
\caption{Variation of stress intensity factors $K_I$ and $K_{II}$ with crack angle, $\beta$.}
\label{fig:incl_crk_Sifs}
\end{figure}

The influence of crack angle $\beta$ on the SIFs is shown in~\fref{fig:incl_crk_Sifs}. A structured mesh $(100 \times 100)$ is used for the study. For crack angles, $0 < \beta < 90$, the crack face intersects the elements in such a way that one region of the split elements (either above or below the crack face) is a polygon. The polygonal subdomain is mapped onto a unit disk (see \fref{fig:split}) instead of subdividing it into triangles. It is seen from~\fref{fig:incl_crk_Sifs} that the numerical results are comparable with the analytical solution.

\subsection{Multiple cracks in tension}
\begin{figure}[htbp]
\centering
\includegraphics[angle=0,width=0.5\textwidth]{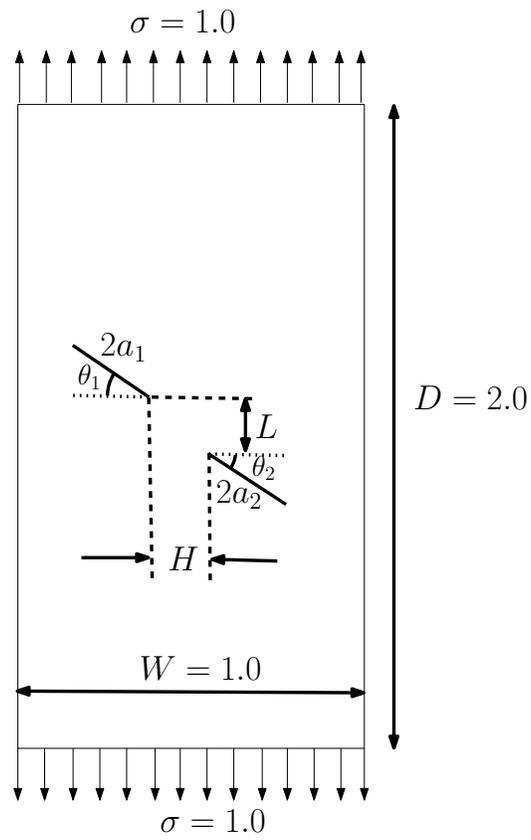}
\caption[]{Plate with multiple cracks: geometry and loads. The length of the cracks are $2a_1$ and $2a_2$. $\theta_1$ and $\theta_2$ are the angles subtended by the crack faces with the horizontal.}
\label{fig:multiCrk}
\end{figure}

In the next example, we consider a plate with two cracks.
The geometry and boundary conditions of the problem are shown
in~\fref{fig:multiCrk}. In this case, the problem is solved only by
the new proposed method. The material properties are: Young's
modulus $E=3\times10^7$ and Poisson's ratio, $\nu=0.3$. A mesh size
of $72 \times 144$ is used for the current study with crack size,
$2a_1=0.2$. The length of the other crack $2a_2$ is varied.

\begin{figure}[htbp]
\centering
\includegraphics[angle=0,width=0.8\textwidth]{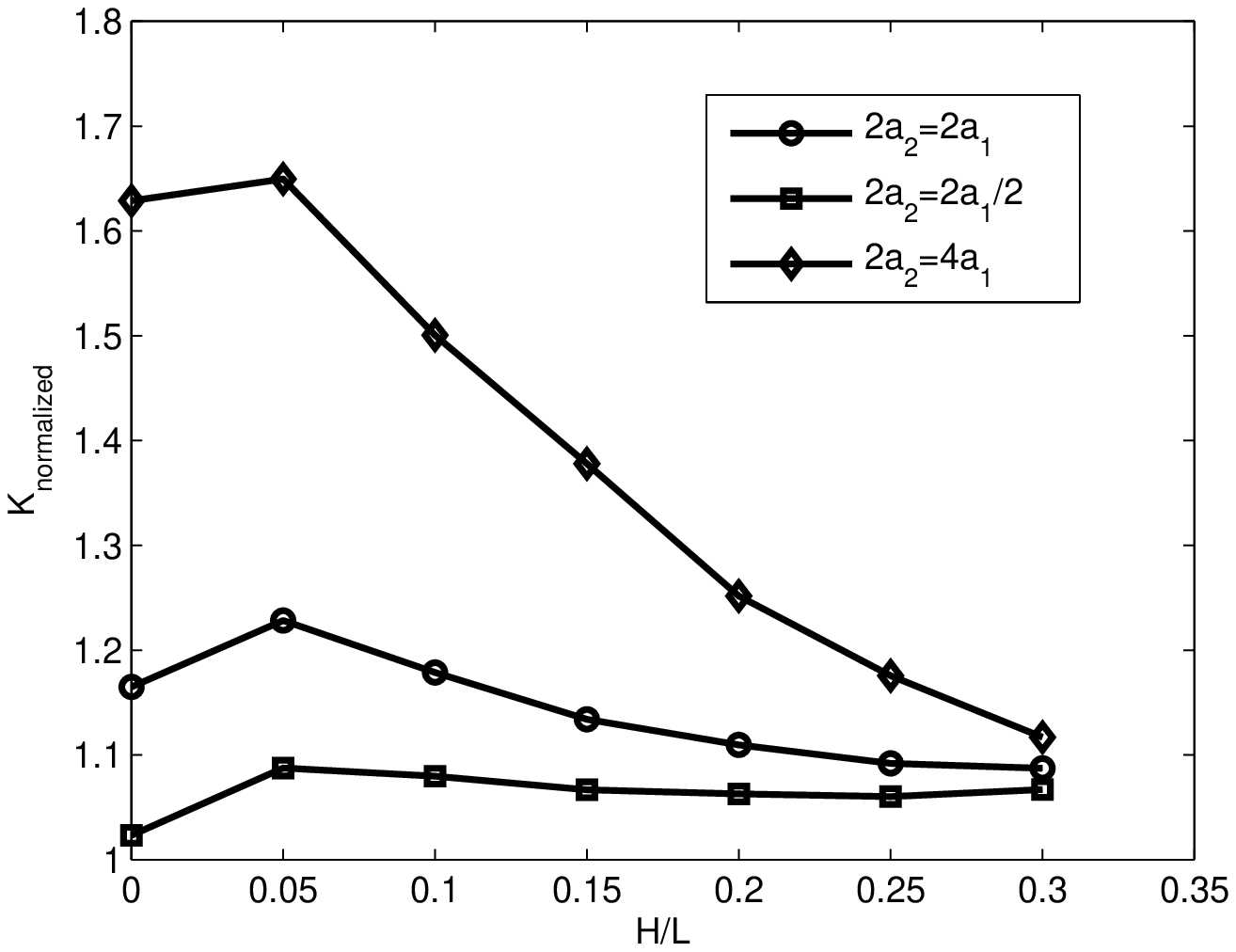}
\caption[]{Influence of $H/L$ ratio and $2a_2/2a_1$ on computed
stress intensity factor.} \label{fig:multicrack1}
\end{figure}

\fref{fig:multicrack1} shows the variation of the mode I SIF for
different $H/L$ ratios and for different ratio of crack
lengths\footnote{crack length ratio = $\frac{2a_2}{2a_1}$} with $\theta_1=0$ and $\theta_2=0$, where $\theta_1$ and $\theta_2$ are the angles subtended by the crack faces with the horizontal (see \fref{fig:multiCrk}). The
normalized mode I SIF\footnote{$K_{{normalized}} =
\frac{K_{\textup{numerical}}}{K_{\textup{analytical}}}$} is plotted
for point $A$ in~\fref{fig:multiCrk}. This is done to
non-dimensionalize the results. Also, the value of $K_\textup{analytical}$
is taken as the value for a plate with a center crack given by

\begin{equation}
K_I = \sigma \sqrt{\pi a \sec\left( \frac{\pi a}{2w} \right)}
\end{equation}

\noindent where $a$ is the half crack length, $w=\frac{W}{2}$ is the
half width of the plate, and $\sigma$ is the far field tensile load
applied at the top of the plate.

It can be seen that with increase in the $H/L$ ratio, the
interaction between the cracks decreases, which is as expected. It is also seen from
~\fref{fig:multicrack1} that as the ratio of crack lengths increases, the normalized mode I SIF also increases. And as the $H/L$ ratio increases, all the curves approach the analytical solution $(K_\textup{analytical})$. The values will not be equal to the analytical solution, because the analytical solution is computed for the case with one center crack.

\begin{figure}[htpb]
\centering
\subfigure[mode I]{\includegraphics[width=0.7\textwidth]{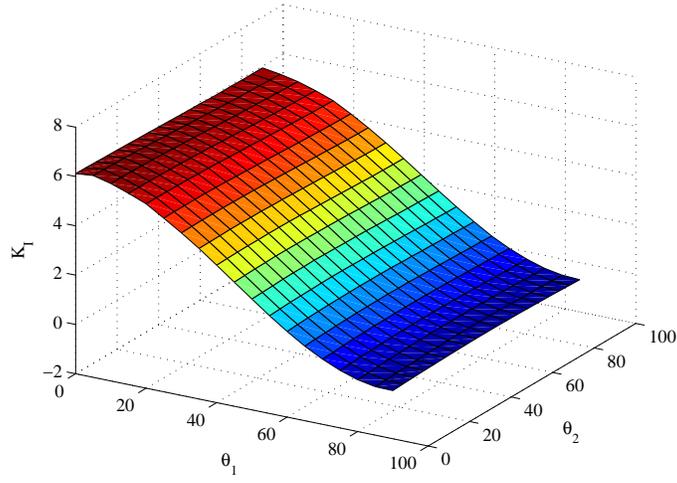}}
\subfigure[mode II]{\includegraphics[width=0.7\textwidth]{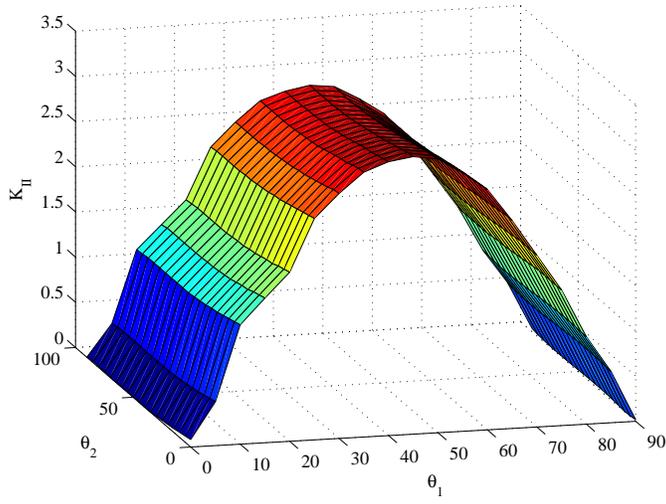}}
\caption{Plate with two cracks: the variation of mode I SIF and mode
II SIF with respect to angle between the cracks for crack lengths
$2a_1=0.2$ and $2a_2=0.2$. The distance between the cracks are:
$H=0.1$ and $L=0.2$} \label{fig:sif_angle}
\end{figure}

Next, the influence of the relative angle between the cracks on the
mode I and the mode II SIF for a crack length of $2a_1=2a_2=0.2$,
with the distance between the cracks: $H=0.1$ and $L=0.2$ is
studied. \fref{fig:sif_angle} shows the variation of mode I and mode
II SIF with the angle subtended by the cracks to the horizontal axis. It
is seen that with the increase in the angle of the crack,
($\theta_1$ and $\theta_2$), the mode I SIF decreases and approaches
zero for $\theta_1=\theta_2=90^o$. While the mode II SIF, initially
increases with increase in the crack angle and reaches the maximum
for the crack angle $\theta_1=\theta_2=45^o$ and decreases with
further increase in the angle.

\subsection{Bimaterial bar}
To illustrate the effectiveness of the proposed method to
integrate weak discontinuities, we simulate a one-dimensional
bimaterial bar discretized with 2D quadrilateral elements. Consider a
two-dimensional square domain $\Omega= \Omega_1 \cup \Omega_2$ of
length $L=2$ with the material interface $\Gamma$ located at
$b=L/2$. The Young's modulus and Poisson's ratio in $\Omega_1 =
(-1,b)\times(-1,1)$ are $E_1=1$, $\nu=0$, and that in $\Omega_2 =
(b,1)\times(-1,1)$ are $E_2=10$, $\nu=0$. With no body forces, the
exact displacement solution with $u_y=0$ at $y=-1$ and $u_y=1$ at
$y=1$ is given by:

\begin{equation}
\renewcommand\arraystretch{1.5}
u(y) = \left\{ \begin{array}{lr} \left(y+1\right)\alpha,&  -1\le y \le b, \\
1+\frac{E_1}{E_2}\left(y-1\right)\alpha,& b\le y \le 1 \end{array}
\right.
\end{equation}

\noindent where,

\begin{equation}
\alpha = \frac{E_2}{E_2(b+1) - E_1(b-1)}
\end{equation}

\begin{figure}[htpb]
\centering
\subfigure[]{\includegraphics[width=0.7\textwidth]{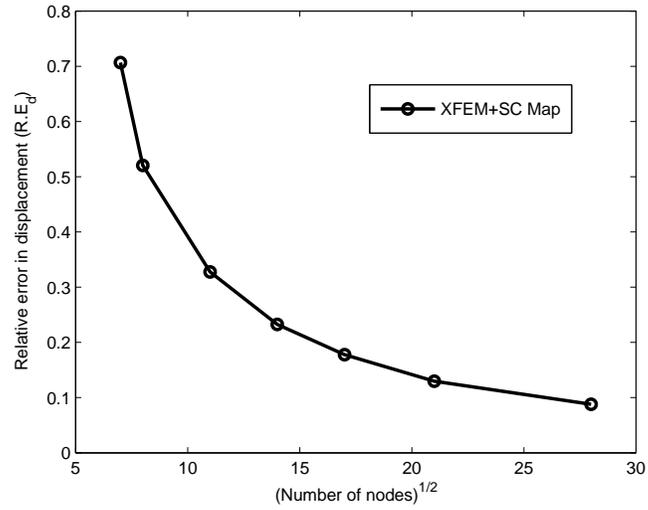}}
\subfigure[]{\includegraphics[width=0.7\textwidth]{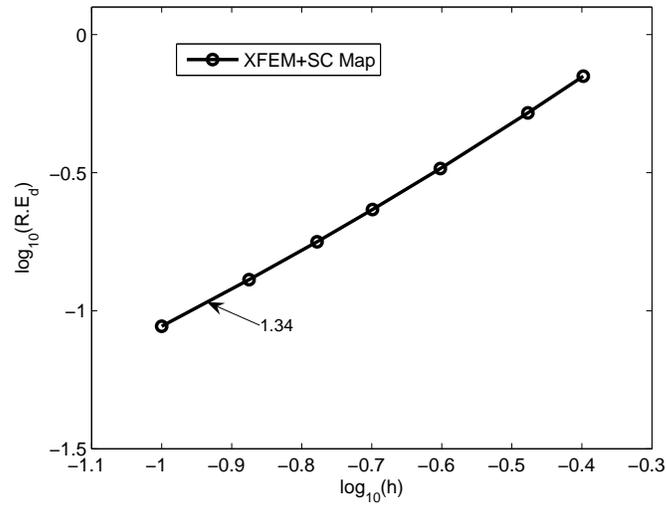}}
\caption{Bimaterial problem: convergence in the displacement
$\left(L_2\right)$ norm: (a) the relative error and (b)
the convergence rate. The method obtains convergence rate of 1.34 in the $L_2-$ norm} \label{fig:bimat_disp}
\end{figure}

The relative error in displacement norm used to measure the accuracy
of the results:

\begin{equation}
R.E_d = \sqrt{\frac{\sum_{i=1}^{\textup{ndof}}
\left(u_i^h-u_i^{\textup{exact}}\right)^2}{\sum_{i=1}^{\textup{ndof}}\left(u_i^{\textup{exact}}\right)^2}}
\times 100 \label{eqn:relerror}
\end{equation}

\fref{fig:bimat_disp} shows the relative error and the rate of
convergence in the displacement norm. It is seen that with decrease
in mesh density, the relative error in displacement norm also
decreases. The rate of convergence is also shown in
the~\fref{fig:bimat_disp}. The proposed method obtains convergence rate of 1.34 in the $L_2-$ norm, very close to the value reported in the literature~\cite{Fries2008}. The convergence rate is suboptimal due to the absence of treatment of the blending elements~\cite{chessablending,Xiao5,B'echet2005,Laborde2005}.

\begin{figure}[htpb]
\centering
\subfigure[]{\includegraphics[angle=0,width=0.4\textwidth]{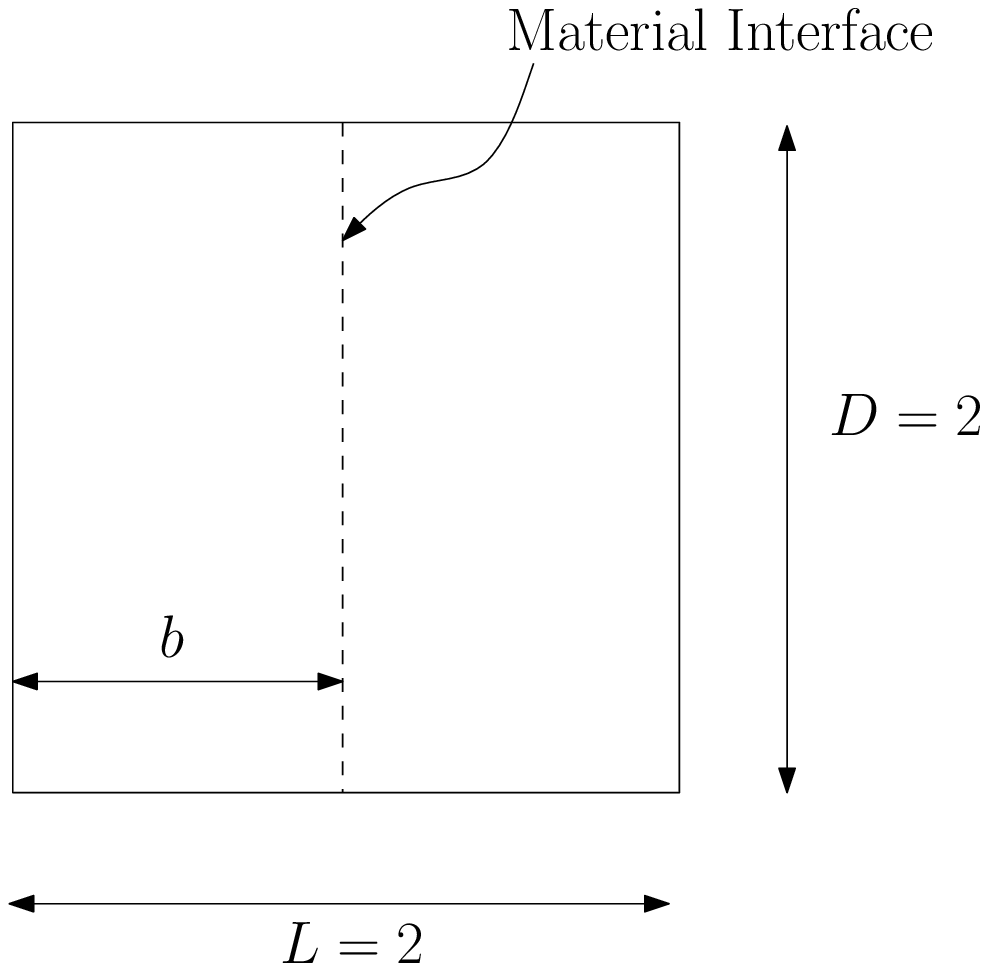}}
\subfigure[]{\includegraphics[angle=0,width=0.4\textwidth]{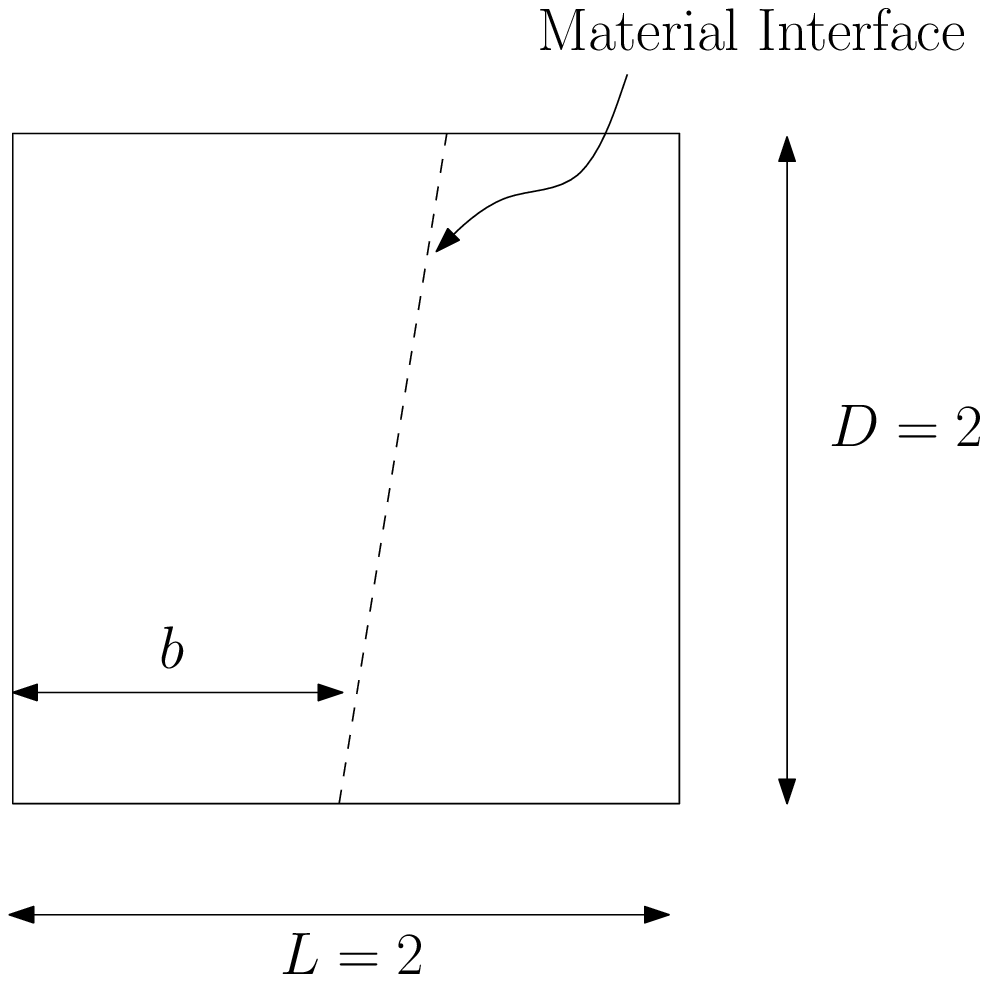}}
\subfigure[]{\includegraphics[angle=0,width=0.4\textwidth]{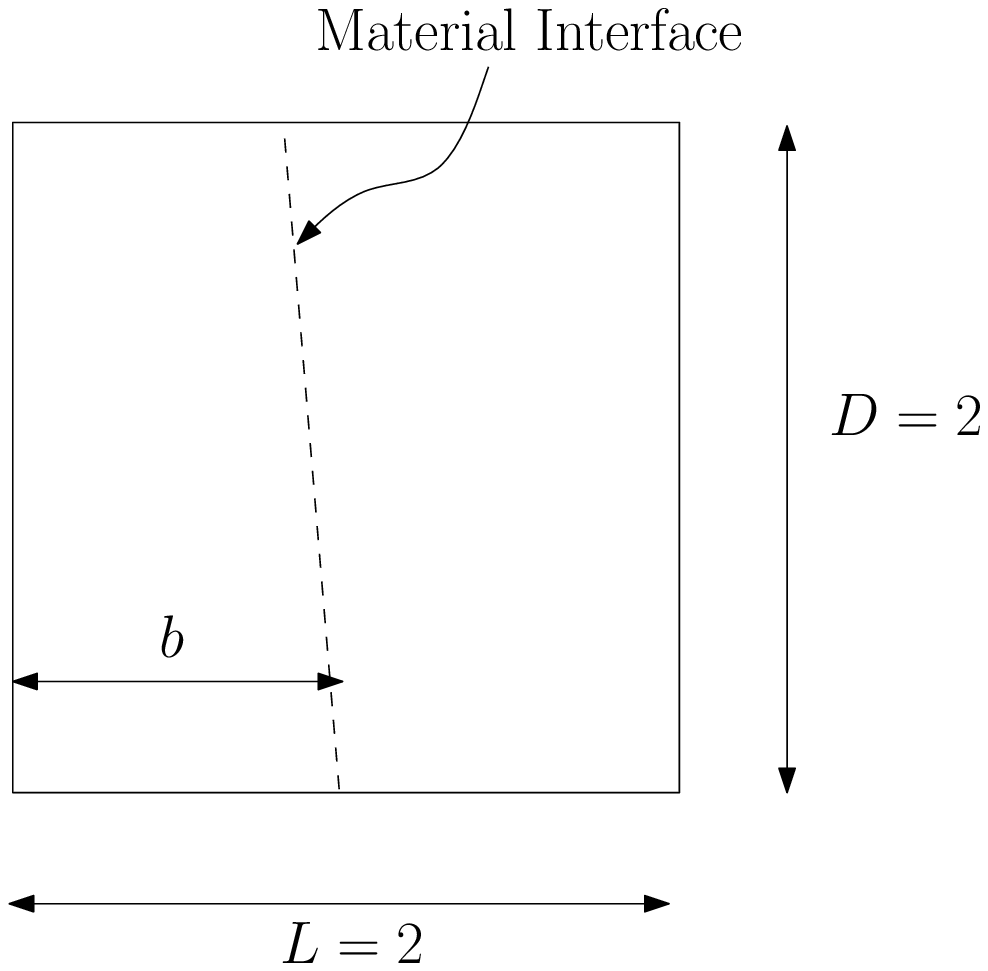}}
\caption{Bimaterial problem: (a)Straight Interface; (b) Slanted Interface (Positive slope) and (c) Slanted Interface (Negative slope)}
\label{fig:matInt}
\end{figure}

Next, the influence of the material interface is studied. Three different configurations of the material interface are considered for the study (see \fref{fig:matInt}).

\begin{figure}[htpb]
\centering
\includegraphics[angle=0,width=0.7\textwidth]{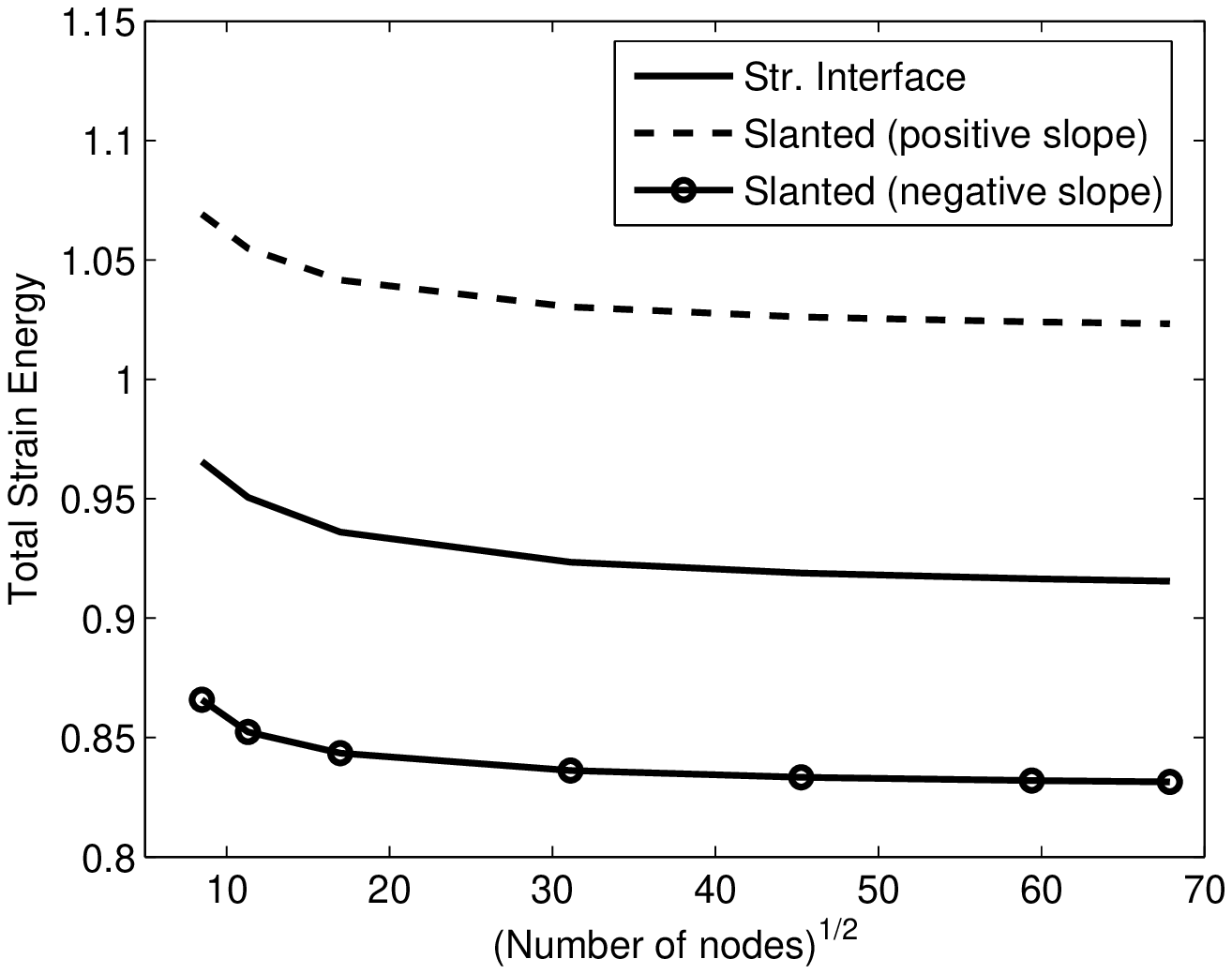}
\caption{Bimaterail problem: convergence of the strain energy with mesh refinement}
\label{fig:bimat_strain}
\end{figure}

\fref{fig:bimat_strain} shows the strain energy convergence with mesh refinement for three different configurations. It is seen that the strain energy converges as the mesh is refined. 

\subsection{Double cantilever beam}
\begin{figure}[htpb]
\centering
\includegraphics[width=0.7\textwidth]{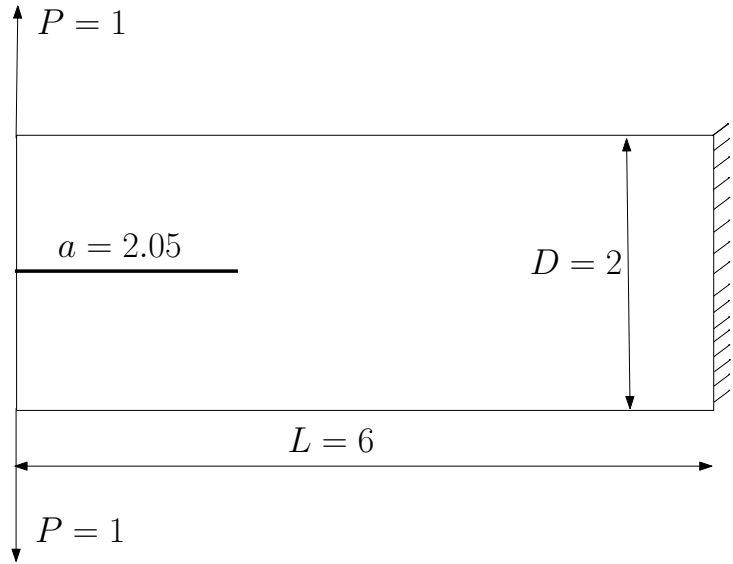}
\caption{Geometry and loads of a double cantilever beam}
\label{fig:doublecb}
\end{figure}

The dimensions of the double cantilever beam (see \fref{fig:doublecb}) are $6\times 2$ and the initial pre-crack with length of $a=2.05$ is considered. The material properties are taken to be Young's modulus, $E=100$ and Poisson's ratio, $\nu=0.3$. And the load $P$ is taken to be unity. By symmetry, a crack on the mid-plane of the beam is under pure mode I and the crack would propagate in a straight line, however, due to small perturbations in the crack geometry, the crack takes a curvilinear path~\cite{Belytschko1999}. A quasi-static crack growth is considered in this study and the growth is governed by the maximum hoop stress criterion~\cite{Erdogan1963}, which states that the crack will propagate from its tip in the direction $\theta_c$ where the circumferential (hoop) stress $\sigma_{\theta \theta}$ is maximum. The critical angle is computed by solving the following equation:

\begin{equation}
K_I \sin(\theta_c) + K_{II} (3\cos(\theta_c)-1) = 0
\label{eqn:prop1}
\end{equation}

\noindent Solving \Eref{eqn:prop1} gives the crack propagation angle~\cite{Sukumar2003}

\begin{equation}
\theta_c = 2 \arctan \left[ { -2\left( {K_{II} \over K_I} \right) \over 1+ \sqrt{1+8 \left({K_{II} \over K_I} \right)^2} } \right]
\end{equation}

\begin{figure}[htpb]
\centering
\subfigure{\includegraphics[width=0.7\textwidth]{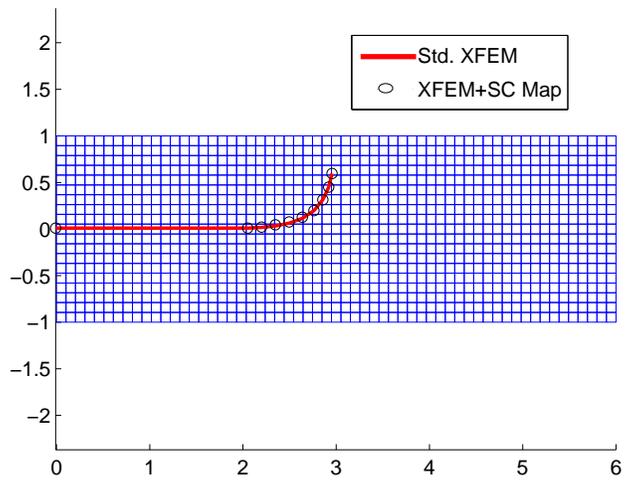}}
\subfigure{\includegraphics[width=0.7\textwidth]{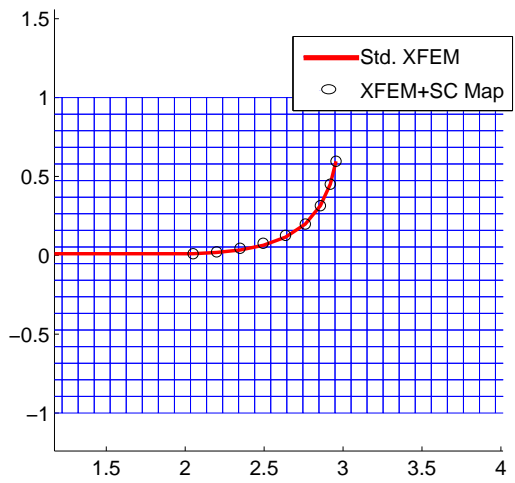}}
\caption{Double cantilever beam: comparison of crack path between the two numerical integration methods.}
\label{fig:crkPath}
\end{figure}

\noindent The crack growth increment, $\Delta a$ is taken to be 0.15 for this study and the crack growth is simulated for 8 steps. The domain is discretized with a structured mesh consisting of 1200 elements. The crack path is simulated using both methods and is shown in~\fref{fig:crkPath}. The crack path qualitatively agrees with the published results~\cite{Belytschko1999}.

\section{CONCLUSION}
In this paper, we used the new numerical integration
proposed for arbitrary polygons in~\cite{Natarajan2009} to integrate
the discontinuous and singular integrands appearing in the XFEM
stiffness matrix. The proposed method eliminates the need to
sub-divide elements cut by strong or weak discontinuities or
containing the crack tip. With a few examples from linear elastic
fracture mechanics and a bimaterial problem, the effectiveness of
the proposed method is illustrated. It is seen that for similar
number of integration points, the proposed technique slightly
outperforms the conventional integration method based on
sub-division. With mesh refinement, both integration techniques
provide convergence of the SIFs to the analytical SIFs. It seems
possible that the proposed technique could serve as a way to
integrate discontinuous approximations in the context of 3D problems
as well, which will be the topic of forthcoming communications.

{\bf Acknowledgement}
The first author acknowledges the financial support of (1) the Overseas Research Students Awards Scheme; (2) the Faculty of Engineering, for period Jan. 2009 - Sept. 2009 and of (3) the School of Engineering (Cardiff University) for the period Sept. 2009 onwards.

The last author gratefully acknowledge the financial support of the Royal Academy of Engineering and of the Leverhulme Trust for Senior Research Fellowship (2009-2010) [Towards the Next Generation Surgical Simulators] http://www.raeng.ork.uk/research/reseacher/leverhulme/current.htm.

\bibliographystyle{wileyj}
\bibliography{mesh_free_reference,sfem_reference,fem_reference,xfem_reference}

\end{document}